\documentclass[11pt, reqno]{amsart}
\usepackage{amsmath}
\usepackage[utf8]{inputenc}
\usepackage{amsmath, amsthm, amsfonts, amssymb, color}
\usepackage{setspace}
\usepackage{mathrsfs}
\usepackage{amsfonts, amsmath}
\usepackage{amsmath,amstext,amsthm,amssymb,amsxtra}
\usepackage{newunicodechar}
\newunicodechar{，}{,}
\newunicodechar{−}{-}
\usepackage{hyperref}
\hypersetup{colorlinks=true,linkcolor=black, urlcolor=blue, citecolor=blue}

\usepackage{color}

\newcommand{\EQ}[1]{\begin{equation}\begin{split} #1 \end{split}\end{equation}}
\newcommand{\EQN}[1]{\begin{equation*}\begin{split} #1 \end{split}\end{equation*}}
\newcommand{\CAS}[1]{\begin{cases} #1 \end{cases}}

\newcommand{\case}[1]{
\EQ{
\begin{array}{l}
\left\{
\begin{aligned}
#1
\end{aligned}
\right.
\end{array}
}
}

\usepackage[margin=3cm, a4paper]{geometry}

\definecolor{deepgreen}{cmyk}{1,0,1,0.5}

\newcommand {\norm}[1]{\left\|#1\right\|}

\newtheorem{theorem}{Theorem}[section]
\newtheorem{proposition}[theorem]{Proposition}

\newtheorem{lemma}[theorem]{Lemma}

\newtheorem{conjecture}[theorem]{Conjecture}
\newtheorem{problem}[theorem]{Problem}
\theoremstyle{remark}
\newtheorem{remark}[theorem]{Remark}

\theoremstyle{definition}

\allowdisplaybreaks

\def\ep{\varepsilon}

\def\th2{\frac{\theta}{2}}

\def\normo#1{\left\|#1\right\|} 
\def\normb#1{\big\|#1\big\|}

\def\norm#1{\|#1\|}

\def\wh#1{\widehat{#1}}

\newcommand{\Div}{{\rm div\,}}

\newcommand{\R}{{\mathbb R}}

\newcommand{\Z}{{\mathbb Z}}
\newcommand{\X}{{\mathbb X}}

\newcommand{\les}{{\lesssim}}

\newcommand{\supp}{{\mbox{supp}}}

\newcommand{\hP}{{\dot{P}}} 

\numberwithin{equation}{section}

\begin{document}

\title[Compressible Navier-Stokes equation]
{
Heat-wave semigroup estimates and the compressible Navier-Stokes Equations
}

\author[R. Chen]{Rongshan Chen}
\address{School of Mathematics, Monash University, Melbourne, VIC 3800, Australia.}
\email{Rongshan.Chen@monash.edu}
\author[Z. Guo]{Zihua Guo}
\address{School of Mathematics, Monash University, Melbourne, VIC 3800, Australia.}
\email{Zihua.Guo@monash.edu}

\thispagestyle{empty}

\begin{abstract}
We study the connection between the heat-wave semigroup and the compressible Navier-Stokes (CNS) equations. On one hand, the parabolic-hyperbolic nature of the CNS motivated the study of space-time estimates for the heat-wave semigroup by using recent harmonic analysis development, in particular, the local smoothing estimates for the wave equations; on the other hand, we apply the new estimates to obtain global well-posedness for CNS in the optimal critical Besov spaces. For the latter purpose, we introduce a new iteration framework for CNS that exploits some intrinsic structures of the equations, so that we can successfully apply these estimates.
\end{abstract}


\maketitle

\tableofcontents

\section{Introduction}

Many PDE models in fluid dynamics exhibit complicated and mixed properties of parabolic, hyperbolic or dispersive types at different frequency scales. The study of each type has seen flourish developments in the last several decades and many adaptive tools were established. For PDE of the mixed types, it is crucial to exploit the combinations of these tools.  That usually involves delicate studies of the mutual influences of different tools.  The purpose of this paper is to demonstrate a framework that connects the heat-wave semigroup estimates and the compressible Navier-Stokes equation. We believe it may work for other parabolic-dispersive PDE. 

\subsection{Compressible Navier-Stokes equations}

We consider the Cauchy problem to the barotropic compressible Navier-Stokes equations on $\mathbb{R}^d$:
\begin{equation}\label{eq:CNS}
\left\{
\begin{aligned}
& \partial_{t}\rho+{\rm div}\, (\rho u)=0,\quad  t>0,\,  x\in \mathbb{R}^{d},\\
&  \partial_{t}(\rho u)+{\rm div}\, (\rho u\otimes u)-\mathcal{A}u+\nabla P(\rho)=0,\quad  t>0,\,  x\in \mathbb{R}^{d},\\
& \rho(0,x)=\rho_{0}(x),\,   u(0,x)=u_{0}(x),\quad x\in \mathbb{R}^{d},  \\
\end{aligned}
\right.
\end{equation}
where $\rho(t,x)$ and $u(t,x)$ represent the fluid's density and velocity, $P(\rho)$ represents the pressure of the fluid  depending on the density and $\mathcal{A} = \mu\Delta + (\mu + \lambda)\nabla \mathrm{div}$ is the Lamé operator, representing the viscosity. The first equation is derived from the conservation of mass and the second equation is the momentum equation with viscosity. 
The study of the Cauchy problem \eqref{eq:CNS} is of fundamental importance in both mathematics and physics. We refer to \cite{Nash,MN1980, PL1998, BaChDa11} and references therein.

In order to simplify the presentations and highlight our methods, we consider only the most canonical and simplest setting. We assume $\mu$ and $\lambda$ are constants and satisfy the elliptic conditions: $\mu > 0$ and $2\mu +\lambda > 0$, and $\gamma = P'(1)=2\mu +\lambda=1$.  We also assume that the fluid density $\rho$ is a small perturbation near $1$ so that $\rho$ has a lower bound. Let $a=\rho-1$ and define
\EQ{
R_0(a) = \frac{1}{a+1},\quad R_1(s) =& \frac{1}{s^2}\int_0^s \frac{P'(1+\tau)-1-\tau}{(1+\tau)}d\tau.
}
Then we can rewrite \eqref{eq:CNS} as:
\begin{equation}\label{eq:CNS2}
\left\{
\begin{aligned}
& \partial_ta+\Div u=f(a,u),\\
& \partial_{t}u-\mathcal{A}u+\nabla a=g(a,u),
\\
& a(0,x)=a_{0}(x),\,   u(0,x)=u_{0}(x)  \\
\end{aligned}
\right.
\end{equation}
where $f,g$ are the nonlinear terms
\case{\label{eq:nonlinearityfg}
f(a,u)=&-\Div(au),\\
g(a,u)=&-u\cdot \nabla u-R_0(a) a\mathcal{A}u-\nabla (R_1(a)a^2).
}

Consider the linear compressible Navier-Stokes equation
\begin{equation}\label{eq:linearCNS}
\left\{
\begin{aligned}
& \partial_t a+\Div u=f,\\
&  \partial_{t}u-\mathcal{A}u+ \nabla a=g.
\end{aligned}
\right.
\end{equation}
We have (cf. \cite{HoZu1995,FD2010,CMZ2010})
\begin{align}\label{linear solution}
\begin{pmatrix} a \\ u \end{pmatrix} 
&= G(t)\begin{pmatrix} a_{0} \\ u_{0} \end{pmatrix} 
+ \int_{0}^{t} G(t-s)\begin{pmatrix} f \\ g \end{pmatrix}  \mathrm{d}s,
\end{align}
where $G(t)$ is the Fourier multiplier operator defined by the symbol 
\begin{equation}\label{Gtsymbol}
\widehat{G}(\xi,t) \triangleq 
\begin{pmatrix}
\dfrac{\lambda_{+}e^{\lambda_{-}t} - \lambda_{-}e^{\lambda_{+}t}}{\lambda_{+}-\lambda_{-}} 
& i\dfrac{e^{\lambda_{-}t} - e^{\lambda_{+}t}}{\lambda_{+}-\lambda_{-}}\xi^\top \\[2ex]
-i\dfrac{e^{\lambda_{-}t} - e^{\lambda_{+}t}}{\lambda_{+}-\lambda_{-}}\xi 
& \dfrac{\lambda_{+}e^{\lambda_{+}t} - \lambda_{-}e^{\lambda_{-}t}}{\lambda_{+}-\lambda_{-}}\dfrac{\xi\xi^\top}{|\xi|^2}
\end{pmatrix}.
\end{equation}
Here $\lambda_\pm$ are the eigenvalues
\EQ{
\lambda_{\pm}=\lambda_\pm(|\xi|) = -\frac{1}{2}|\xi|^2 \pm \frac{1}{2}\sqrt{|\xi|^4 - 4|\xi|^2}.
}
One can see that, for high-frequency $|\xi|>2$, the semigroup $e^{t\lambda_\pm(D)}$ is parabolic, and 
\EQ{
\lambda_+\sim -|\xi|^2,\quad \lambda_-\sim -1, \quad \mbox{as}\quad |\xi|\to \infty.
}
For low frequency $|\xi|<2$, the semigroup $e^{t\lambda_\pm(D)}$ is parabolic-dispersive,
\EQ{
\lambda_\pm(|\xi|)=-\frac{|\xi|^2}{2}(1\pm iS(|\xi|)), \qquad S(|\xi|)=\sqrt{\frac{4}{|\xi|^2}-1}. 
}
Note that in the transition frequency region $||\xi|-2|\ll 1$, $G(t)$ is not singular and behaves essentially like $e^{Ct\Delta}$.

\subsection{Heat-wave semigroup estimates}

It is natural to see that the estimates for the parabolic-dispersive semigroup $e^{t\lambda_\pm(D)}$ would play a crucial role in the study of the equation \eqref{eq:CNS2}.
It was found by the second named author-Song-Yang (\cite{GSY2025}) that when acting on low frequency
\EQ{
e^{t\lambda_\pm (D)}\approx e^{ct\Delta+iCtD}
}
which is the heat-wave semigroup, by proving: for $1\leq p\leq \infty$ and $k\leq 0$
\EQ{\label{l1estimate} 
e^{-Ct 2^{2k}}\norm{e^{\mp 2it D}\hP_k f}_{L_x^p}\les \|e^{t\lambda_\pm (D)}\hP_{k}f\|_{L_x^{p}}\lesssim  e^{-ct 2^{2k}}\norm{e^{\mp 2it D}\hP_k f}_{L_x^p}.
}
Here $\hP_k$ is the Littlewood-Paley dyadic projection operator.
The heat semigroup and wave operator are the most canonical objects in parabolic PDE, and hyperbolic/dispersive PDE. They were extensively studied individually. 
However, these separate results are not trivially applied to the $e^{t\lambda_\pm (D)}$ as one can imagine that there are mutual influences, competition and balance between the dissipative and dispersive behaviours. In this paper we study systematically the space-time estimates for the heat-wave semigroup, in particular, for our purposes, the following estimates: for $k\leq 0$
\EQ{\label{eq:L2Lq}
\|e^{t\Delta+itD}\hP_{k}f\|_{L_t^{2}L_x^{q}(\R\times \R^d)}
\les& 2^{k\delta(q,r,d)}\norm{\hP_kf}_{L_x^{r}}.
}
Note that $e^{t\Delta+itD}\hP_{k}\approx e^{-ct2^{2k}} e^{itD}\hP_k$. So \eqref{eq:L2Lq} is equivalent to the time-weighted wave operator estimates:
\EQ{\label{eq:L2Lq2}
\|e^{-ct2^{2k}}e^{itD}\hP_{k}f\|_{L_t^{2}L_x^{q}(\R\times \R^d)}
\les& 2^{k\delta(q,r,d)}\norm{\hP_kf}_{L_x^{r}}
}
which is very delicate. 

By using the fixed-time estimates for the wave operator $e^{itD}$ due to Peral \cite{Peral} (Section IV) and Miyachi \cite{Miyachi}: for $1<p<\infty$
\EQ{\label{eq:Waveoperator}
\norm{e^{itD} f}_{L^p(\R^d)}\leq C_t \norm{f}_{H^{s_p}_p(\R^d)}
}
with $s_p=(d-1)|\frac{1}{2}-\frac{1}{p}|$, it was proved in \cite{GSY2025} that: for $1\leq p\leq \infty, k\leq 0$
\EQ{\label{eq:semiLqest}
\|e^{-ct2^{2k}}e^{itD}\hP_{k}f\|_{L^{p}(\R^d)}\lesssim e^{-ct 2^{2k}}2^{-ks_p}\|\hP_{k}f\|_{L^{p}(\R^d)},
}
and then integrating in time, one can get: for $q\in [1,\infty], k\leq 0$
\EQ{\label{eq:semiLqest2}
\|e^{t\Delta+itD}\hP_{k}f\|_{L_t^qL_x^{p}(\R\times \R^d)}\lesssim 2^{-2k/q}2^{-ks_p}\|\hP_{k}f\|_{L^{p}(\R^d)}.
}
The required regularity $s_p$ for \eqref{eq:Waveoperator} is sharp, which is the main reason for the loss of regularity in \eqref{eq:semiLqest2}. Thus, for $q=\infty$, \eqref{eq:semiLqest2} is sharp.

In 1991, Sogge (\cite{S1991}) observed smoothing effect if taking integration over time, and formulated the famous 
{\bf local smoothing conjecture}: 
we have
\EQ{\label{eq:localsmoothingest}
\norm{e^{itD} f}_{L_{t\in [1,2]}^pL_x^p(\R^d)}\leq C \norm{f}_{H^{s_p-\sigma}_p(\R^d)}
}
for some $\sigma$. See Conjecture \ref{Conjecture3.1}. This conjecture is a major open problem in harmonic analysis. One can see that \eqref{eq:L2Lq2} is closely related to \eqref{eq:localsmoothingest} as $e^{-ct2^{2k}}$ plays the role as temporal localization. Another important class of space-time estimates is the Strichartz estimates, namely \eqref{eq:L2Lq2} when $r=2$.
In Section 3, we will use the local smoothing estimates and Strichartz estimates to study \eqref{eq:L2Lq2}. The new estimates improve the loss of regularity for \eqref{eq:semiLqest2} when $q<\infty$. The improvement is crucial for our purposes to study equation \eqref{eq:CNS2}.

\subsection{Main results and strategy of proof}

We study the well-posedness of \eqref{eq:CNS2} in the critical spaces. 
Note that the following scaling transform
\begin{equation*}
\begin{aligned}
(a_{0}, u_{0})\rightsquigarrow(a_{0}(\tau x),\tau u_{0}(\tau x)),\quad(a(t,x), u(t,x))
\rightsquigarrow(a(\tau^{2}t, \tau x),\tau u(\tau^{2}t, \tau x)), \quad \tau>0
\end{aligned}
\end{equation*}
leaves equation \eqref{eq:CNS2} invariant provided that the pressure law  $P$ is replaced by $\tau^{2}P$. 
A natural choice of space which is scaling invariant for \eqref{eq:CNS2} is the critical Besov space
\begin{equation}\label{criticalspace}
(a_{0},\, u_{0})\in \X_p:=\dot{B}_{p,1}^{\frac{d}{p}} \times \dot{B}_{p,1}^{-1+\frac{d}{p}}.
\end{equation}
Under the incompressible limit (see \cite{RH2015}), \eqref{eq:CNS} formally reduces to the incompressible Navier-Stokes equation:
\case{\label{eq:inns}
\partial_t u - \Delta u  + \nabla  P =& -u\cdot \nabla u, \\
\nabla \cdot u  =& 0, \\
u (0,x) =& u_{0}
}
The equation \eqref{eq:inns} is invariant under the scaling transform $u(t,x) \to \lambda u(\lambda^{2}t, \lambda x)$. The critical Besov space for \eqref{eq:inns} is $\dot{B}_{p,1}^{-1+{d/p}}$.

For the equation \eqref{eq:inns} which is parabolic, global well-posedness with small data in $\dot B^{-1+d/p}_{p,q}$ was proved in \cite{Cannone1997, P1996, C1992} for $1 \leq p < \infty, 1 \leq q \leq\infty$ by perturbative methods.  Ill-posedness in $\dot B^{-1}_{\infty,q}$ in the sense that the solution map is discontinuous was shown in \cite{B2008} for $q=\infty$, and in \cite{T2010} for $q>2$ and then in \cite{BW2015} for $1\leq q\leq 2$.
It is natural and interesting to extend these results to the compressible Navier-Stokes equation \eqref{eq:CNS2}.
However, \eqref{eq:CNS2} is a coupled system of hyperbolic and parabolic equations. The perturbative methods may not apply and the results turn out to be very delicate. 
For initial data in $\X_p$, local existence for $1\leq p<2d$ and uniqueness for $1\leq p\leq d$ was established in \cite{D2001,D2005} (See also \cite{BaChDa11}). The uniqueness for $d<p<2d$ was proved in \cite{D2014} via a Lagrangian approach for \eqref{eq:CNS2}. 
On the other hand, the range $p<2d$ is optimal for well-posedness in $\X_p$. Ill-posedness in the sense that the continuity of the solution map $S_T$ fails at the origin in $\X_p$ was proved in \cite{CMZ2015} for $p>2d$, and then in \cite{IO2022} for $p=2d$. In \cite{GYZ2024}, the second named author-Yang-Zhang proved the continuity of the solution map $S_T$ from $\X_p \to C_T\X_p$ for $1\leq p<2d$ via a combination of the Langrangian approach and the method of frequency envelope, and thus the local well-posedness is completely settled.

For global results, in a seminal paper \cite{D2000}, by delicate energy method, Danchin proved global existence and uniqueness for small initial data 
$(a_0,u_0)\in ( \dot{B}_{2,1}^{d/2-1} \cap \dot{B}_{2,1}^{d/2}) \times \dot{B}_{2,1}^{d/2-1}$. Additional conditions on the low frequency seem necessary in view of properties of the linear equations \eqref{eq:linearCNS}. Later, the results were extended to $L^p$-framework in \cite{FD2010, CMZ2010, H20112}, proving global existence and uniqueness with 
initial data satisfying
\EQ{
\|(a_{0},u_{0})^\ell\|_{\dot{B}_{2,1}^{d/2-1}} 
+ \|(a_{0},u_0)\|_{\X_p} \ll 1,
}
where \(p\) satisfies $2 \leq p \leq \min\left\{4, 
\frac{2d}{d-2}\right\}$ for $d \geq 3$ or $2 \leq p < 4$ for $d = 2$.
Note that uniqueness was restricted to $2\leq p\leq d$ in \cite{FD2010,CMZ2010}, and was improved in \cite{H20112} by introducing an effective velocity. 
Here for function $f$, $f^\ell$ (and $f^h$) means the low-frequency (high-frequency) component. 
Compared to the local results (LWP holds for $p<2d$), the range of $p$ for global results is much smaller. 
Recently, for $d=3$, the second-named author, Song and Yang \cite{GSY2025} proved global well-posedness for all $p<6$ assuming
\EQ{
\|(a_{0},u_{0})^\ell\|_{\dot{B}_{3,1}^{-2/3}} 
+ \|(a_{0},u_0)\|_{\X_p} \ll 1.
}
The key ingredients of \cite{GSY2025} are: 1) $L^q$-framework for the low-frequency and \eqref{eq:semiLqest2}; 2) A nonlinear transform $J$ applied to \eqref{eq:CNS2}: 
\EQ{\label{eq:Jtran}
(a,u) \rightarrow J(a,u) :=& (a,\tilde{u}),\\
\tilde{u} :=& m^{\ell}+u^{h} = u+ (au)^{\ell}.
}
The $J$-transform was shown to be a bijection in the resolution spaces. It is problematic for the low frequency if the regularity is too negative. In particular, the method of \cite{GSY2025} is not sufficient for obtaining optimal well-posedness in dimensions four and higher. 

In this paper we continue to study the global well-posedness in higher dimensions, by using the new heat-wave semigroup estimates. We introduce a new iteration framework which can obtain global well-posedness for \eqref{eq:CNS2} in $\X_p$ with the optimal range of $p$ in the dimensions $d\leq 4$ and wider range for $d\geq 5$. It seems to us that this is a natural formulation for \eqref{eq:CNS2} and much more flexible than the one used in \cite{GSY2025}. It may be useful for other purposes. Now we discuss our strategy. 

A major obstacle for \eqref{eq:CNS2} is caused by the negative regularity for the low frequency. For example, to deal with the difference estimate for the local well-posedness, the $J$-transform and the estimate for parabolic-dispersive semigroup developed in \cite{GSY2025}, one of the main constraints is from handling the negative regularity near $0$-frequency. Unlike for the high-frequency, one does not have tools like smoothing effect for low frequency component.
Inspired by \cite{GSY2025,HoZu1995,HoZu1997}, we consider the equivalent momentum equation $(a,m)$ for \eqref{eq:CNS}: $m=(1+a)u$
\begin{equation}\label{eq:CNSm}
\left\{
\begin{aligned}
& \partial_{t}a+\operatorname{div}m=0,\quad  t>0,\,  x\in \mathbb{R}^d,\\
&  \partial_{t}m-\mathcal{A}m+\nabla a=h(a,m),\quad  t>0,\,  x\in \mathbb{R}^d,\\
& \rho(0,x)=\rho_{0}(x),\,   m(0,x)=m_{0}(x),\quad x\in \mathbb{R}^d.  \\
\end{aligned}
\right.
\end{equation}
where
$h(a,m)$ denotes the nonlinear terms $h(a,m)=\sum_{i=1}^{3}h_i(a,m)$
with
\case{\label{eq:nonlinearh123}
h_1(a,m)
&=-\operatorname{div}\big(R_0(a)m\otimes m\big),\\
h_2(a,m)
&=-\mathcal{A}\big(R_0(a)a m\big),\\
h_3(a,m)
&=\nabla\big(R_2(a)a^2\big).
}
Where $R_2(a)=\frac{P(1+a)-P(1)-a}{a^2}$.
Note that compared to the velocity equation \eqref{eq:CNS2}, the nonlinearity $h$ is of some divergence forms which improves the low-frequency estimate for $m$. The nonlinearity $h_2(a,m)$ contains a term $R_0(a)m \mathcal{A}a$ which is problematic to the high-frequency estimate for $m$. It is natural to use the effective velocity method for the high frequency. Thus as the $J$-transform for equation \eqref{eq:CNS2}, we may introduce a nonlinear transform $K$ for equation \eqref{eq:CNSm}:
\EQ{\label{eq:Ktran}
(a,m) \rightarrow K(a,m) :=& (a,\tilde{m}),\\
\tilde{m} :=& m^{\ell}+u^{h} = m+ (\frac{-a}{a+1}m)^{h}.
}
Namely, we still use momentum for the low-frequency and velocity for the high-frequency. However, the $K$-transform is problematic for the high-frequency estimates, as the highest regularity needed for $m$ is $\dot B_{p,1}^{d/p+1}$ while the regularity for $a$ is $\dot B_{p,1}^{d/p}$. This means $\tilde m$ can not stay in the same spaces as $m$, and then we couldn't get closed estimates using the momentum equation \eqref{eq:CNSm}.

Another important way to deal with the obstacle is to reduce the loss of regularity in the heat-wave semigroup estimates. In \cite{GSY2025}, \eqref{eq:semiLqest2} was used. \eqref{eq:semiLqest2} is sharp when $q=\infty$, however, there is an improvement when $q<\infty$. It is natural to use \eqref{eq:L2Lq2}, namely $L_t^2L^p$ type estimate, which is weaker than \eqref{eq:semiLqest2}. 
Our new framework is as follows: instead of using a single equation, we use both the velocity and momentum equations and the two nonlinear transforms $J,K$ which can combine the strength of the two formulations. More precisely, we consider the following system for $(a^\ell,m^\ell, a^h,u^h,u^\ell,m^h)$
\begin{equation}\label{eq:CNSmu}
\left\{
\begin{aligned}
&\partial_t a^\ell+\operatorname{div}m^\ell=0,\\
&\partial_t m^\ell-\mathcal{A}m^\ell+\nabla a^\ell=[h(a,m)]^\ell,\\
&\partial_t a^h+\operatorname{div}u^h=[f(a,u)]^h,\\
&\partial_t u^h-\mathcal{A}u^h+\nabla a^h=[g(a,u)]^h\\
&u^\ell=m^\ell-(au)^\ell\\
&m^h=u^h+(au)^h\\
\end{aligned}
\right.
\end{equation}
By this formulation, we have much more flexibilities in designing resolution spaces. In particular, we do not need to require that $J$ or $K$ transform is invertible in the full resolution spaces. Instead, we can use much weaker spaces for $u^\ell$ than $m^\ell$, and for $m^h$ than $u^h$.  For example, our new formulation \eqref{eq:CNSmu} allows us to use $L_t^\infty \dot B^{d/p}_{p,1}$ for $u^\ell$, instead of $L_t^\infty \dot B^{d/q-1}_{q,1}$. Then we can successfully apply \eqref{eq:L2Lq2} and close all the estimates. This is the key difference and our main novelty compared to \cite{GSY2025}.

Our main result is:

\begin{theorem}\label{thm}
Let $d\geq 2$. Assume 
$2\leq q\leq p\leq 2q,\, q<\frac{3d}{d-1}, \, p<2d$. Assume $0 < \delta \ll 1$ and $(a_0, u_0) \in E_\delta$ where (note that $m=(1+a)u$)
\begin{equation}
E_\delta = \{(a_0, u_0) \in \mathbb{X}_p : \|(a_0,m_0)^\ell\|_{\dot{B}^{\frac{3d}{2q}-\frac{1+d}{2}-}_{q,1}} + \|(a_0,u_0)\|_{\mathbb{X}_p}\leq \delta\}.
\end{equation}
Then the Cauchy problem \eqref{eq:CNS} has a unique global solution $(a,u)$ such that $(a, m,u) \in X_{q,p}$ satisfying $\|(a, m,u)\|_{X_{q,p}} \leq C\delta$ where $X_{q,p}$ is the Banach space defined by the norm: 
\EQ{\label{eq1.22}
&\|(a, m,u)\|_{X_{q,p}} \\
=& \|(a, m)^\ell\|_{\widetilde{L}^\infty_t(\dot{B}^{\frac{d}{q}}_{q,1}) \cap \widetilde{L}^2_t(\dot{B}^{\theta}_{q,1})} +\|u^\ell\|_{\widetilde{L}^\infty_t(\dot{B}^{\frac{d}{q}}_{q,1}) \cap \widetilde{L}^2_t(\dot{B}^{\frac{d}{q}}_{q,1})}\\ &+\|a^h\|_{\widetilde{L}^\infty_t(\dot{B}^{\frac{d}{p}}_{p,1}) \cap \widetilde{L}^1_t(\dot{B}^{\frac{d}{p}}_{p,1})}+ \|u^h\|_{\widetilde{L}^\infty_t(\dot{B}^{-1+\frac{d}{p}}_{p,1}) \cap \widetilde{L}^1_t(\dot{B}^{1+\frac{d}{p}}_{p,1})}+\|m^h\|_{\widetilde{L}^\infty_t(\dot{B}^{\frac{d}{p}-1}_{p,1})\cap \widetilde{L}^2_t(\dot{B}^{\frac{d}{p}}_{p,1})}
}
where $\theta=\frac{d}{2q}$.
\end{theorem}

\begin{remark}\label{rem:explain}
We can take $q=p/2$, then the above theorem holds for $p<\min(\frac{6d}{d-1},2d)$. It covers the optimal range for high frequency $p<2d$ for $2\leq d\leq 4$, and improves the range obtained in the previous results (\cite{FD2010, CMZ2010, H20112}) for $d\geq 5$. 
\end{remark}

\medskip

\noindent {\bf Notations and definitions.}
We collect some notations, definitions and some nonlinear estimates. Some were used in \cite{GSY2025} and we continue to use in this paper. 

$\bullet$ $A \lesssim B$ means that $A \leq CB$ for some constant $C$. $A\sim B$ means $A\les B$ and $B\les A$.

$\bullet$ For $a\in \R$, $a-$ denotes $a-\epsilon$ for any sufficiently small $\epsilon>0$. Similarly we define $a+$.

$\bullet$ The Fourier transform of $f \in \mathcal{S}$ (the Schwartz class) is defined as
\[
\widehat{f}(\xi) = \mathcal{F}[f](\xi) :=(2\pi)^{-d/2} \int_{\mathbb{R}^{d}} f(x) e^{-i \xi \cdot x}  dx.
\]

$\bullet$ Fix a real-valued $C_0^\infty(\R^n)$ function $\eta(\xi)$ such that: $\eta$ is a radially symmetric and radially decreasing; $\supp \eta\subset \{\xi\in \R^n:|\xi|\leq 2^{0.01}\}$; and $0\leq \eta(\xi)\leq
1$, $\eta(\xi)\equiv 1$ if $|\xi|\leq 1$.

$\bullet$ Let $\psi(\xi)=\eta(\xi/2)-\eta(\xi)$. For $k\in \Z$, we define
$\psi_k(\xi)=\psi(2^{-k}\xi)$, $\eta_0=\eta$, $\eta_k=\psi_k$ for $k\geq 1$, $\eta_k=0$ for $k\leq -1$.

$\bullet$ We define the Littlewood-Paley projectors:
for $k\in \Z$
\begin{align*}
\wh{\dot{P}_kf}(\xi)=\psi_k(\xi)\hat{f},\,\,\wh{{P}_kf}(\xi)=\eta_k(\xi)\hat{f},
\,\, \wh{{P}_{\leq k}f}(\xi)=\eta(2^{-k}\xi)\hat{f}, \,\, P_{\geq k}=I-P_{\leq k-1}.
\end{align*}
By the definition, we see that $\dot P_{\leq j-2}u\dot P_jv$ has Fourier support $\subset \{\xi:2^{j-2}\leq |\xi|\leq 2^{j+2}\}$.

$\bullet$ Throughout this paper, for \( z \in \mathcal{S}'(\mathbb{R}^d) \), we use $z^{\ell}$, $z^h$ to denote the low-frequency and high-frequency component respectively, namely
\begin{equation}\label{hilo}
z^{\ell}:=P_{\leq j_{0}-1}z \quad \text{and} \quad 
z^{h}:=z-z^{\ell},
\end{equation}
where $j_0$ is a universal parameter to be determined later. $(a,u)^{\ell}$ (resp $(a,u)^h$) means $(a^{\ell},u^{\ell})$ (resp $(a^h,u^h)$).

$\bullet$ The norm of homogeneous Besov space: for $1\leq p,r \leq \infty$
\begin{equation}\label{eq2.4}
\begin{aligned}
& \norm{f}_{\dot{B}_{p,r}^{s}} := \norm{2^{js}\norm{\hP_{j}f}_{L^{p}(\mathbb{R}^d)}}_{l^{r}(\mathbb{Z})}.
\end{aligned}
\end{equation}

$\bullet$ We denote the function spaces ${L^{q}_{T}(\dot{B}_{p,r}^{s})}$ and ${\tilde{L}^{q}_{T}(\dot{B}_{p,r}^{s})}$ to be Banach spaces with the following norms:
\EQN{
\|u(t,\cdot)\|_{L^{q}_{T}(\dot{B}_{p,r}^{s})}=&\|\|u(t,\cdot)\|_{\dot{B}_{p,r}^{s}}
\|_{L^{q}_{T}},\\
\|u(t,\cdot)\|_{\tilde{L}^{q}_{T}(\dot{B}_{p,r}^{s})} =& \big\| \big( 2^{js} \|\hP_{j}u(t,\cdot)\|_{L_{T}^{q}(L_x^{p})} \big) \big\|_{\ell^{r}(\mathbb{Z})}.
}
The index $T$ will be omitted when $T = +\infty$. We denote by $\mathcal{C}_{b}(\dot{B}_{p,r}^{s})$ the subset of functions in ${L_t}^{\infty}(\dot{B}_{p,r}^{s})$ that are also continuous from $\mathbb{R}_{+}$ to $\dot{B}_{p,r}^{s}$. 

\section{Space-time estimates for the heat-wave semigroup}

In this section, we derive the global-in-time estimates for the semigroup $e^{t\Delta+itD}$ acting on the low frequency. First, for the heat semigroup, we have the following classical estimates:
\begin{lemma}[\cite{BaChDa11}]\label{lem:heatsemigroup}
Assume $u$ solves the heat equation:
\EQ{
\partial_t u - \mu\Delta u    =&  f,\\
u(0,x)=& u_0(x)
}
where $\mu>0$ is a constant. Then for $s\in\mathbb{R}$ and $\rho_{1}, \rho,p, r\in [1,\infty]$ with $\rho_{1}\geq \rho$, we have
\EQ{\label{eq:heatsemigroup}
\mu^{1/\rho_{1}}\|u\|_{\tilde{L}_{T}^{\rho_{1}}\dot{B}_{p,r}^{s+\frac{2}{\rho_{1}}}}
\les\|u_{0}\|_{\dot{B}_{p,r}^{s}}
+\mu^{1/\rho-1}\|f\|_{\tilde{L}_{T}^{\rho}\dot{B}_{p,r}^{s-2+\frac{2}{\rho}}}.
}
\end{lemma}
The estimate \eqref{eq:heatsemigroup} is scaling invariant. The wave operator $e^{itD}$ has different scaling than $e^{t\Delta}$, and thus affects the estimates. We recall that by using fixed-time estimates for the wave operator \eqref{eq:Waveoperator}, the following estimates were obtained in \cite{GSY2025}:
for $p,q\in [1,\infty]$
\EQ{\label{eq:semiLqest3}
\|e^{t\Delta+itD}\hP_{k}f\|_{L_t^qL_x^{p}(\R\times \R^d)}\lesssim 2^{-2k/q}2^{-ks_p}\|\hP_{k}f\|_{L^{p}(\R^d)}
}
where
\EQ{
s_p=(d-1)\big|\tfrac{1}{2}-\tfrac{1}{p}\big|.
}
Inspired by these results, we propose to study the following problem:
\begin{problem}
Find the optimal conditions for $(r,q,p,s)$ such that the estimate
\EQ{\label{eq:heat-wave}
\bigl\|e^{t\Delta+itD}\hP_{k}f\bigr\|_{L_t^q L_x^p}
\les 2^{ks}\bigl\|\hP_k f\bigr\|_{L_x^r}
}
holds for all $f$ and $k\leq 0$.
\end{problem}
The estimate \eqref{eq:heat-wave} is equivalent to 
\EQ{\label{eq:heat-wave2}
\bigl\|e^{-t2^{2k}}e^{itD}\hP_{k}f\bigr\|_{L_t^q L_x^p}
\les 2^{ks}\bigl\|\hP_k f\bigr\|_{L_x^r}
}
which is closely related to the local smoothing conjecture proposed by Sogge (\cite{S1991}):
\begin{conjecture}[local smoothing conjecture]\label{Conjecture3.1}
Assume $d\geq 2$. Then for $\ep>0$ we have
\EQ{\label{eq:localsmoothing}
\norm{e^{itD} f}_{L_{t\in [1,2]}^pL_x^p(\R^d)}\les \norm{f}_{H^{s+\ep}_p(\R^d)}
}
if $s=\gamma_p$ with
\EQ{
\gamma_p=
\CAS{
s_p-\frac{1}{p}=\tfrac{d-1}{2}-\tfrac{d}{p}, \quad &\frac{2d}{d-1}\leq p<\infty;\\
0, &\quad 2<p< \frac{2d}{d-1}.
}
}
\end{conjecture}

When $d=2$, the conjecture was recently proved positively by \cite{GWZ}. When $p\geq \frac{2(d+1)}{d-1}$, the conjecture was verified, while the case $\frac{2d}{d-1}\leq p<\frac{2(d+1)}{d-1}$ remains open for $d\geq 3$. \footnote{Proof of the conjecture was released by OpenAI for $d=3$ (24 Sep 2026), for all dimensions (26 Oct 2026). Our results for $d\geq 4$ do not rely on these results, but on Lemma \ref{lem:localsmoothingest4}.}

\begin{lemma}[Local smoothing estimates, \cite{BG2015}]\label{lem:localsmoothingest4}
Assume $d\geq 2$. Then for $\ep>0$ we have \eqref{eq:localsmoothing} holds for $s=\beta_p$ where
\EQ{
\beta_p=
\CAS{
s_p-\frac{1}{p}=\tfrac{d-1}{2}-\tfrac{d}{p}, \quad &\frac{2(d+1)}{d-1}\leq p<\infty;\\
\tfrac{s_p}{2}=\tfrac{d-1}{2}(\tfrac{1}{2}-\tfrac{1}{p}), &\quad 2<p< \frac{2(d+1)}{d-1}.
}
}
\end{lemma}

We use the above lemma to get improved semigroup estimates than \eqref{eq:semiLqest3} in the case $q\leq p$. The other case can also be improved using interpolation.  

\begin{lemma}[Improved semigroup estimates]
Assume $k\leq 0$, $1\leq r\leq q\in [1,\infty)$ and $q\geq 2$. Then we have
\begin{equation}\label{eq:semi-qq}
\bigl\|e^{t\lambda_\pm (D)}\hP_{k}f\bigr\|_{L_t^r L_x^q}
\les 2^{-k(\tfrac{2}{r}+\beta_q+)}\bigl\|\hP_k f\bigr\|_{L_x^q}.
\end{equation}
\end{lemma}
\begin{proof}

Given a function $f$, for $j\in \Z$, we denote $f_{j}(x)=f(x/2^{j})$.
By the properties of Fourier transform, we have
\begin{equation}\label{95}
\begin{aligned}
[e^{-th(2^j D)}\hP_{k}f](x)=&\mathcal{F}^{-1}(e^{-th(2^j\xi)}\hat{f}(\xi)\chi(\frac{\xi}{2^{k}}))\\
=&\mathcal{F}^{-1}(e^{-th(2^{j}\xi)}\hat{f}(2^{-j}{2^j\xi})
\chi(\frac{2^j\xi}{2^{j+k}}))\\
=&2^{-jn}[\mathcal{F}^{-1}(e^{-th(\xi)}\hat{f}(2^{-j}\xi)
\chi(\frac{\xi}{2^{j+k}}))](2^{-j}x)\\
=&[e^{-th(D)}\hP_{j+k}f_{-j}](2^{-j}x).\\
\end{aligned}
\end{equation}  
By \eqref{l1estimate} we get for $1\leq q\leq \infty$ and $k\leq 0$
\EQN{
\|e^{t\lambda_\pm (D)}\hP_{k}f\|_{L_x^{q}}\lesssim  e^{-ct 2^{2k}}\norm{e^{\mp 2it D}\hP_k f}_{L_x^q}.
}
Then 
\EQ{
\|e^{t\lambda_\pm (D)}\hP_{k}f\|_{L_t^{r}L_x^{q}}\lesssim&  \normo{\normb{e^{-ct 2^{2k}}\norm{e^{\mp 2it D}\hP_k f}_{L_x^q}}_{L_{t\in [2^j,2^{j+1}]}^r}}_{l_{j\in \Z}^r}\\
\les& \normo{e^{-\tilde c  2^{j+2k}}\normb{\norm{e^{\mp 2it D}\hP_k f}_{L_x^q}}_{L_{t\in [2^j,2^{j+1}]}^r}}_{l_{j\in \Z}^r}\\
\les& \normo{e^{-\tilde c  2^{j+2k}}2^{j/r}\norm{e^{\mp 2it2^j D}\hP_k f}_{L_{t\in [1,2]}^rL_x^q}}_{l_{j\in \Z}^r}\\
\les& \normo{e^{-\tilde c  2^{j+2k}}2^{j/r}2^{dj/q}\norm{e^{\mp 2it D}\hP_{j+k} f_{-j}}_{L_{t\in [1,2]}^rL_x^q}}_{l_{j\in \Z}^r}.
}
This reduces to the estimate $\norm{e^{\mp 2it D}\hP_{j+k} f_{-j}}_{L_{t\in [1,2]}^rL_x^q}$. We will apply the local smoothing estimates.

{\bf Case 1.} $j+k\leq 0$

By the H\"older inequality in time and Taylor expansion, we get
\EQ{
\norm{e^{\mp 2it D}\hP_{j+k} f_{-j}}_{L_{t\in [1,2]}^rL_x^q}\les 2^{-jd/q}\norm{f}_{L^q}.
}

{\bf Case 2.} $j+k\geq 1$

By the H\"older inequality in time and applying Lemma \ref{lem:localsmoothingest4}, we have
\EQ{
\norm{e^{\mp 2it D}\hP_{j+k} f_{-j}}_{L_{t\in [1,2]}^rL_x^q}\les& \norm{e^{\mp 2it D}\hP_{j+k} f_{-j}}_{L_{t\in [1,2]}^qL_x^q}\\
\les& 2^{(j+k)^+(\beta_q+)}\norm{f_{-j}}_{L_x^q}\\
\les& 2^{(j+k)^+(\beta_q+)}2^{-jd/q}\norm{f}_{L_x^q}.
}
Therefore, we get
\EQ{
\|e^{t\lambda_\pm (D)}\hP_{k}f\|_{L_t^{r}L_x^{q}}
\les& \normo{e^{-\tilde c  2^{j+2k}}2^{j/r}2^{(j+k)^+(\beta_q+)}}_{l_{j\in \Z}^r}\norm{f}_{L_x^q}\\
\les& 2^{-k(\tfrac{2}{r}+\beta_q+)}\norm{f}_{L_x^q}.
}
\end{proof}

Another important class of estimates is the Strichartz estimates (see \cite{Guo2018}). 

\begin{lemma}[Strichartz estimates]\label{lem:str}
Assume $d\geq 2$, $2\leq r, q\leq \infty$ satisfies the {\bf wave-admissible} conditions: 
\EQ{
\tfrac{1}{r}\leq \tfrac{d-1}{2}(\tfrac{1}{2}-\tfrac{1}{q}), \quad (r,q,d)\neq (2,\infty,3).
}
Then we have
\EQN{
\|e^{it\sqrt{-\Delta}}f\|_{L^{r}_{t}(\R;L^{q}_{x})}\les \|f\|_{\dot{H}^{s}}
}
where $s=\alpha(r,q):=\frac{d}{2}-\frac{d}{q}-\frac{1}{r}$. 
\end{lemma}

By the above Lemma and \eqref{eq:heat-wave}, we get the following estimates: for wave-admissible pair $(r,q)$, $k\in \Z$
\EQ{\label{eq:Str}
\|e^{t\Delta+itD}\hP_{k}f\|_{L_t^{r}L_x^{q}}
\les& 2^{k\alpha(r,q)}\norm{\hP_kf}_{L_x^2}.
}
Then by interpolation between \eqref{eq:Str} and \eqref{eq:semi-qq}, we can obtain a family of estimates \eqref{eq:semiLqest3}.

\section{Priori estimates for general dimension}

In this section, we derive the global-in-time a-priori estimates for $(a,u,m)$ using the following system
\begin{equation}\label{eq3.1}
\left\{
\begin{aligned}
&\partial_t a^\ell+\operatorname{div}m^\ell=0,\\
&\partial_t m^\ell-\mathcal{A}m^\ell+\nabla a^\ell=[h(a,m)]^\ell,\\
&\partial_t a^h+\operatorname{div}u^h=[f(a,u)]^h,\\
&\partial_t u^h-\mathcal{A}u^h+\nabla a^h=[g(a,u)]^h\\
&u^\ell=m^\ell-(au)^\ell\\
&m^h=u^h+(au)^h\\
\end{aligned}
\right.
\end{equation}
where $f,g,h$ are given by \eqref{eq:nonlinearityfg} and \eqref{eq:nonlinearh123}.
The resolution spaces for $a,m,u$ are Banach spaces defined with the following norms:
\EQ{
\norm{a}_A&=\norm{a^\ell}_{\widetilde L_t^\infty(\dot B_{q,1}^{d/q})\cap\widetilde L_t^2(\dot B_{q,1}^{\theta})}
+\|a^h\|_{\widetilde L_t^\infty(\dot B_{p,1}^{\frac dp})\cap \widetilde L_t^1(\dot B_{p,1}^{\frac dp})}:=A^\ell+A^h,\\
\norm{m}_M&=\norm{m^\ell}_{\widetilde L_t^\infty(\dot B_{q,1}^{d/q})\cap\widetilde L_t^2(\dot B_{q,1}^{\theta})}
+\|m^h\|_{\widetilde L_t^2(\dot B_{p,1}^{\frac{d}{p}})}:=M^\ell+M^h,\\
\norm{u}_U&=\|u^\ell\|_{\widetilde L_t^\infty(\dot B_{p,1}^{\frac{d}{p}})\cap \widetilde L_t^2(\dot B_{p,1}^{\frac{d}{p}})}
+\|u^h\|_{\widetilde L_t^\infty(\dot B_{p,1}^{-1+\frac dp})\cap\widetilde L_t^1(\dot B_{p,1}^{1+\frac dp})}:=U^\ell+U^h.
}
We also denote $A=\norm{a}_A,M=\norm{m}_M,U=\norm{u}_U$, $\theta$ is a parameter to be determined.

We will use bootstrap arguments. We assume throughout this section the following condition:
\EQ{
A+M+U\leq \epsilon_0\ll 1.
}
We assume the estimates: for $d\geq 1$ and $k\leq 0$, $q\geq 2$
\EQ{\label{eq:L2Lq}
\|e^{t\lambda_\pm (D)}\hP_{k}f\|_{L_t^{2}L_x^{q}(\R^d)}
\les& 2^{k\delta(q,r,d)}\norm{\hP_kf}_{L_x^{r}}
}
and derive general estimates under conditions depending on $\delta(q,r,d)$ and $\theta$.
The estimates \eqref{eq:L2Lq} were studied in Section 2 and we apply these results to prove Theorem \ref{thm} in the next section . 

\subsection{Estimates for $(a^\ell,m^\ell)$}

In this subsection, we prove the estimates for $(a^\ell,m^\ell)$ using the first two equations in \eqref{eq3.1}. 
By \eqref{eq:L2Lq} we have the following linear estimates.

\begin{lemma} \label{lem:lowlinear}
Assume that $j_{0}>0, \, s \in \mathbb{R}$, $d \geq 1$, $p \in [1, \infty]$. Let $z = (z_{1}, z_{2})$ satisfy 
\begin{equation}
\left\{
\begin{aligned}
& \partial_{t}z_{1}+\operatorname{div} z_{2}=f\\
&  \partial_{t}z_{2}-\mathcal{A} z_{2}+\nabla  z_{1}=g, \\
& (z_{1},z_{2})\big|_{t=0}=(z_{1,0},z_{2,0}). \\
\end{aligned}
\right.
\end{equation}
Then for all $\theta \in \R$, 
\EQN{
\|(z_1,z_2)^\ell\|_{\widetilde L_t^\infty(\dot B_{q,1}^{\frac{d}{q}})\cap \widetilde L_t^2
(\dot B_{q,1}^{\theta})}\les& 
\|(z_{1,0},z_{2,0})^\ell\|_{\dot B_{q, 1}^{\frac{2d-1}{q}+\frac{1-d}{2}}\cap \dot B_{q,1}^{\theta+\delta(q,q,d)})}\\
&+\inf_{r\leq q}\|(f,g)^\ell\|_{\widetilde L_t^1(\dot B_{q,1}^{\frac{d-1}{q}+\frac{1-d}{2}+\frac{d}{r}})}+\inf_{r\leq q}\|(f,g)^\ell\|_{\widetilde L_t^1(\dot B_{r,1}^{\theta+\delta(q,r,d)})}.
}
\end{lemma}
\begin{proof}
The proof is similar to Lemma 3.5 of \cite{GSY2025}.
\end{proof}

We have the following

\begin{proposition}\label{prop:low}
Assume
\EQ{\label{eq:hestcond}
\frac{2d-1}{q}+\frac{3-d}{2}>0, \quad 2\le q\leq p\leq 2q,\quad -\delta(q,q,d)-1<\theta\leq \frac{d}{2q}.
}
Then we have
\EQ{
A^\ell+M^\ell\;\les\; \|(a_0,m_0)^\ell\|_{\dot B_{q,1}^{\frac{2d-1}{q}+\frac{1-d}{2}}\cap \dot B_{q,1}^{\theta+\delta(q,q,d)})}+(A+M)^2.
}
\end{proposition}

By applying Lemma \ref{lem:lowlinear} to the first two equations in the system (\ref{eq3.1}), we can get the following estimates:
\begin{equation}\label{eq3.8}
\begin{aligned}
&
\|(a,m)^\ell\|_{\widetilde L_t^\infty(\dot B_{q,1}^{\frac{d}{q}})}
+
\|(a,m)^\ell\|_{\widetilde L_t^2
(\dot B_{q,1}^{\theta})}
\\
\lesssim&
\|(a_0,m_0)^\ell\|_{\dot B_{q,1}^{\frac{2d-1}{q}+\frac{1-d}{2}}}
+\|(a_0,m_0)^\ell\|_{\dot B_{q,1}^{\theta+\delta(q,q,d)})}\\
&+\inf_{r\leq q}\|h(a,m)^\ell\|_{\widetilde L_t^1(\dot B_{r,1}^{\frac{d-1}{q}+\frac{1-d}{2}+\frac{d}{r}})}+\inf_{r\leq q}\|h(a,m)^\ell\|_{\widetilde L_t^1(\dot B_{r,1}^{\theta+\delta(q,r,d)})}.
\end{aligned}
\end{equation}
It remains to prove the nonlinear estimates.

We will use the following product estimates frequently in the nonlinear estimates.
\begin{lemma}\label{lem:product}
We have\\
\begin{itemize}
\item[(1)] Assume $1\leq p \leq \infty$. Then
\EQ{\label{eq:prodest1}
\norm{fg}_{\dot{B}_{p,1}^{d/p}} \les& \norm{f}_{\dot{B}_{p,1}^{d/p}} \norm{g}_{\dot{B}_{p,1}^{d/p}}.
}
\item[(2)] Assume $2\leq p<2d$. Then
\EQ{\label{eq:prodest2}
\norm{fg}_{\dot{B}_{p,1}^{d/p-1}} \les& \norm{f}_{\dot{B}_{p,1}^{d/p}} \norm{g}_{\dot{B}_{p,1}^{d/p-1}}.
}
\end{itemize}
\end{lemma}
\begin{proof}
The first estimate is the classical algebraic property. For the second estimate, it follows from Bony's paraproduct decomposition and Bernstein's inequality. The restriction $p<2d$ is from the high-high to low interactions. 
\end{proof}

\begin{proposition}\label{prop3.3}
Assume \eqref{eq:hestcond}.
Then
\EQ{
\inf_{r\leq q}\norm{(h(a,m))^{\ell}}_{\widetilde L_t^1
(\dot B_{r,1}^{\frac{d-1}{q}+\frac{1-d}{2}+\frac{d}{r}})\cap \widetilde L_t^1
(\dot B_{r,1}^{\theta+\delta(q,r,d)})}\les (A+M)^2.
}
\end{proposition}

\begin{proof}

First we estimate $\inf_{r\leq q}\norm{(h_{1}(a,m))^{\ell}}_{\widetilde L_t^1
(\dot B_{r,1}^{\frac{d-1}{q}+\frac{1-d}{2}+\frac{d}{r}})\cap \widetilde L_t^1
(\dot B_{r,1}^{\theta+\delta(q,r,d)})}$.
Since $h_{1}(a,m) = \operatorname{div}(R_0(a) m \otimes m)$ and we get 
\EQ{
&\inf_{r\leq q}\norm{(h_{1}(a,m))^{\ell}}_{\tilde{L}_{t}^{1}(\dot B_{r,1}^{\frac{d-1}{q}+\frac{1-d}{2}+\frac{d}{r}})\cap \widetilde L_t^1(\dot B_{r,1}^{\theta+\delta(q,r,d)})}\\
\les& \inf_{r\leq q}\norm{R_0(a)m \otimes m))^{\ell}}_{\tilde{L}_{t}^{1}(\dot B_{r,1}^{\frac{d-1}{q}+\frac{3-d}{2}+\frac{d}{r}})\cap\tilde{L}_{t}^{1}(\dot B_{r,1}^{\theta+\delta(q,r,d)+1})}.
}
We decompose $m\otimes m=m^\ell\otimes m^\ell+2m^\ell\otimes m^h+m^h\otimes m^h$.
There are three types of nonlinear estimates.

For the term $(R_0(a) m^\ell\otimes m^\ell)^\ell$, we take $r=q$. Since $\norm{R_0(a)}_{L_t^\infty L_x^\infty}\les 1$, we get
\EQ{\label{eq:prodest}
&\norm{(R_0(a) m^\ell\otimes m^\ell)^{\ell}}_{\tilde{L}_{t}^{1}(\dot B_{q,1}^{\frac{d-1}{q}+\frac{3-d}{2}+\frac{d}{q}})\cap \tilde{L}_{t}^{1}(\dot B_{q,1}^{\theta+\delta(q,q,d)+1})}\\
\les& \norm{(R_0(a) m^\ell\otimes m^\ell)}_{\tilde{L}_{t}^{1}(\dot B_{q,\infty}^{0})}\\
\leq &\sum_{k_2,k_3\le j_0}\norm{R_0(a)\hP_{k_2} m^\ell \otimes \hP_{k_3} m^\ell}_{{L}_t^1L^{q}}\\
\les &\sum_{k_2,k_3\le j_0} \min(2^{k_2d/q}, 2^{k_3d/q})\norm{\hP_{k_2} m^\ell}_{L_t^2L^q} \norm{\hP_{k_3} m^\ell}_{{L}_t^2L^{q}}\\
\les & \norm{m^\ell}^2_{\tilde{L}_t^2(\dot{B}_{q,1}^{\theta})}\les (M^\ell)^2
}
provided that
\EQ{
\frac{2d-1}{q}+\frac{3-d}{2}>0\\
-\delta(q,q,d)-1<\theta \leq \frac{d}{2q}.\label{eq3.12}
}

Similarly, for the term $(R_0(a) m^h\otimes m^h)^\ell$, we take $r=q$ and get
\EQ{
&\norm{( R_0(a) m^h \otimes m^h)^{\ell}}_{\tilde{L}_{t}^{1}(\dot B_{q,1}^{\frac{d-1}{q}+\frac{3-d}{2}+\frac{d}{q}})\cap \tilde{L}_{t}^{1}(\dot B_{q,1}^{\theta+\delta(q,q,d)+1})}\\
\les& \sum_{k_2,k_3>j_0} \norm{R_0(a)\hP_{k_2} m^h\otimes \hP_{k_3} m^h}_{{L}_t^1L^{q}}\\
\les& \sum_{k_2,k_3>j_0} \norm{\hP_{k_2} m^h}_{L_t^2L^{2q}} \norm{\hP_{k_3} m^h}_{{L}_t^2L^{2q}}
\les (M^h)^2
}
provided that
\EQ{
\frac{2d-1}{q}+\frac{3-d}{2}>0, \quad 2\leq p\leq 2q.\\
\theta+\delta(q,q,d)+1>0.\label{eq3.14}
}

For the term $(R_0(a) m^\ell\otimes m^h)^\ell$, taking $r=q$, we get
\EQ{
&\norm{(R_0(a)  m^\ell \otimes m^h)^{\ell}}_{\tilde{L}_{t}^{1}(\dot B_{q,1}^{\frac{d-1}{q}+\frac{3-d}{2}+\frac{d}{q}})\cap \tilde{L}_{t}^{1}(\dot B_{q,1}^{\theta+\delta(q,q,d)+1})}\\
\les &\sum_{k_2\le j_0,\;j_0<k_3\le j+3}\norm{R_0(a)\hP_{k_2} m^\ell \otimes \hP_{k_3} m^h}_{{L}_t^1L^{q}}\\
\les &\sum_{k_2\le j_0,\;j_0<k_3\le j+3}\norm{\hP_{k_2} m^\ell}_{L_t^2L^{\frac{qp}{p-q}}} \norm{\hP_{k_3} m^h}_{{L}_t^2L^{p}}\\
\les & \norm{m^\ell}_{\tilde{L}_t^2(\dot{B}_{q,1}^{\theta})}\norm{m^h}_{\tilde{L}_t^2(\dot{B}_{p,1}^{\frac{d}{p}})}\les M^\ell \cdot M^h
}
provided that
\EQ{
\frac{2d-1}{q}+\frac{3-d}{2}>0, \quad q\leq p\\
-\delta(q,q,d)-1<\theta\leq \frac{d}{p}.\label{eq3.16}
}

Gathering all the conditions \eqref{eq3.12}, \eqref{eq3.14} and \eqref{eq3.16}, we get 
\EQ{\label{eq:h1est}
\inf_{r\leq q}\norm{(h_{1}(a,m))^{\ell}}_{\widetilde L_t^1
(\dot B_{r,1}^{\frac{d-1}{q}+\frac{1-d}{2}+\frac{d}{r}})\cap \widetilde L_t^1
(\dot B_{r,1}^{\theta+\delta(q,r,d)})}\les M^2.
}
provided that \eqref{eq:hestcond} holds. 

For $h_{2}(a,m) = -\mathcal{A}(R_0(a) a m)$ and $h_{3}(a,m) = \nabla (R_2(a)a^2)$, since $M\les A$, by \eqref{eq:h1est} we have
\EQ{\label{eq:h23est}
\inf_{r\leq q}\norm{(h_{2}(a,m))^{\ell}}_{\widetilde L_t^1
(\dot B_{r,1}^{\frac{d-1}{q}+\frac{1-d}{2}+\frac{d}{r}})\cap \widetilde L_t^1
(\dot B_{r,1}^{\theta+\delta(q,r,d)})}\les& AM\\
\inf_{r\leq q}\norm{(h_{3}(a,m))^{\ell}}_{\widetilde L_t^1
(\dot B_{r,1}^{\frac{d-1}{q}+\frac{1-d}{2}+\frac{d}{r}})\cap \widetilde L_t^1
(\dot B_{r,1}^{\theta+\delta(q,r,d)})}\les& A^2
}
provided that \eqref{eq:hestcond} holds.  Thus the proposition is proved.
\end{proof}

\subsection{Estimates for $(a^h,u^h)$}

In this subsection, we prove the following estimates using the third and fourth equations in \eqref{eq3.1} and the method of effective velocity (see \cite{H20111,GSY2025}). 
\begin{proposition}\label{prop:high}
Assume $2\le q\le p<2d,\, \theta \leq \frac{d}{q}$. Then we have
\begin{equation}
A^h+U^h
\lesssim
\|u_0^h\|_{\dot B_{p,1}^{-1+\frac d p}}
+\|a_0^h\|_{\dot B_{p,1}^{\frac d p}}+(A+U)^2 .
\end{equation}
\end{proposition}
\begin{proof}
Let
\[
\mathbb P:=I+\nabla(-\Delta)^{-1}\operatorname{div},
\qquad
\mathbb Q:=-\nabla(-\Delta)^{-1}\operatorname{div}
\]
denote the orthogonal projectors onto the spaces of divergence-free and potential vector fields, respectively. 

Applying the projections $\mathbb P$ and $\mathbb Q$ to (\ref{eq:CNS2}), we can obtain
\begin{equation}\label{eq3.20}
\left\{
\begin{aligned}
&\partial_t a+\operatorname{div}\mathbb Qu=f,\\
&\partial_t\mathbb Qu-\Delta\mathbb Qu+\nabla a=\mathbb Qg,\\
&\partial_t\mathbb Pu-\mu\Delta\mathbb Pu=\mathbb Pg.
\end{aligned}
\right.
\end{equation}
First, by Lemma \ref{lem:heatsemigroup} we have
\begin{equation}\label{eq3.21}
\|\mathbb Pu^h\|_
{\widetilde L_t^\infty(\dot B_{p,1}^{-1+\frac d p})}
+
\|\mathbb Pu^h\|_
{\widetilde L_t^1(\dot B_{p,1}^{1+\frac d p})}
\lesssim
\|\mathbb Pu_0^h\|_{\dot B_{p,1}^{-1+\frac d p}}
+
\|\mathbb Pg^h\|_
{\widetilde L_t^1(\dot B_{p,1}^{-1+\frac d p})}.
\end{equation}
Then we define the effective velocity $w:=\nabla(-\Delta)^{-1}a+\mathbb Qu$ used by \cite{H20111}. From (\ref{eq3.20}), we get that
\[
\partial_t w-\Delta w
=
\nabla(-\Delta)^{-1}(f-\operatorname{div}g)
+w-\nabla(-\Delta)^{-1}a.
\]
By Lemma \ref{lem:heatsemigroup} and Bernstein's inequality we obtain
\EQ{
&\|w^h\|_{\widetilde L_t^\infty(\dot B_{p,1}^{-1+\frac d p})\cap \widetilde L_t^1(\dot B_{p,1}^{1+\frac d p})}
\\
&\lesssim
\|w_0^h\|_{\dot B_{p,1}^{-1+\frac d p}}
+\|(f,\operatorname{div}g)^h\|_
{\widetilde L_t^1(\dot B_{p,1}^{\frac d p-2})}
+2^{-2j_0}
\Bigl(
\|w^h\|_{\widetilde L_t^1(\dot B_{p,1}^{1+\frac d p})}
+\|a^h\|_{\widetilde L_t^1(\dot B_{p,1}^{\frac d p})}
\Bigr)
\\
&\lesssim
\|\mathbb Qu_0^h\|_{\dot B_{p,1}^{-1+\frac d p}}
+2^{-2j_0}\|a_0^h\|_{\dot B_{p,1}^{\frac d p}}
+2^{-2j_0}
\Bigl(
\|w^h\|_{\widetilde L_t^1(\dot B_{p,1}^{1+\frac d p})}
+\|a^h\|_{\widetilde L_t^1(\dot B_{p,1}^{\frac d p})}
\Bigr)
\\
&\quad +\|(f,\operatorname{div}g)^h\|_
{\widetilde L_t^1(\dot B_{p,1}^{\frac d p-2})}.
}
Taking $j_0$ sufficiently large, we get
\EQ{
&\|w^h\|_{\widetilde L_t^\infty(\dot B_{p,1}^{-1+\frac d p})\cap \widetilde L_t^1(\dot B_{p,1}^{1+\frac d p})}\\
&\lesssim
\|u_0^h\|_{\dot B_{p,1}^{-1+\frac d p}}
+2^{-2j_0}\|a_0^h\|_{\dot B_{p,1}^{\frac d p}}
+2^{-2j_0}\|a^h\|_{\widetilde L_t^1(\dot B_{p,1}^{\frac d p})}+\|(f,\operatorname{div}g)^h\|_
{\widetilde L_t^1(\dot B_{p,1}^{\frac d p-2})}.\label{eq3.23}
}

From (\ref{eq3.20}), we get
\begin{equation}
\partial_t a+u\cdot\nabla a+a
=
-\operatorname{div}w-a\,\operatorname{div}u.
\end{equation}
Then, by Lemma \ref{lem:heatsemigroup}, we get
\begin{equation}\label{eq3.25}
\begin{aligned}
&\|a^h\|_
{\widetilde L_t^\infty(\dot B_{p,1}^{\frac d p})
\cap
\widetilde L_t^1(\dot B_{p,1}^{\frac d p})}
\\
&\lesssim
\|a_0^h\|_{\dot B_{p,1}^{\frac d p}}
+\|(\operatorname{div}w-a\,\operatorname{div}u)^h\|_
{\widetilde L_t^1(\dot B_{p,1}^{\frac d p})}
+\int_0^t
\|\nabla u\|_{\dot B_{p,1}^{\frac d p}}
\|a\|_{\dot B_{p,1}^{\frac d p}}
\,d\tau
\\
&\lesssim
(1+2^{-2j_0})
\|a_0^h\|_{\dot B_{p,1}^{\frac d p}}
+\|u_0^h\|_{\dot B_{p,1}^{-1+\frac d p}}
+2^{-2j_0}
\|a^h\|_{\widetilde L_t^1(\dot B_{p,1}^{\frac d p})}
+\|(a\cdot\operatorname{div}u)^h\|_
{\widetilde L_t^1(\dot B_{p,1}^{\frac d p})}
\\
&\qquad
+\|((f,\operatorname{div}g))^h\|_
{\widetilde L_t^1(\dot B_{p,1}^{\frac d p-2})}.
\end{aligned}
\end{equation}
Taking $j_0$ sufficiently large, we get
\EQ{
\|a^h\|_
{\widetilde L_t^\infty(\dot B_{p,1}^{\frac d p})
\cap
\widetilde L_t^1(\dot B_{p,1}^{\frac d p})}
\lesssim&
\|a_0^h\|_{\dot B_{p,1}^{\frac d p}}
+\|u_0^h\|_{\dot B_{p,1}^{-1+\frac d p}}\\
&+\|(a\cdot\operatorname{div}u)^h\|_
{\widetilde L_t^1(\dot B_{p,1}^{\frac d p})}
+\|((f,\operatorname{div}g))^h\|_
{\widetilde L_t^1(\dot B_{p,1}^{\frac d p-2})}.
}

Therefore, we obtain
\EQ{
&\|a^h\|_
{\widetilde L_t^\infty(\dot B_{p,1}^{\frac d p})
\cap
\widetilde L_t^1(\dot B_{p,1}^{\frac d p})}+\|u^h\|_{\widetilde L_t^\infty(\dot B_{p,1}^{-1+\frac dp})\cap\widetilde L_t^1(\dot B_{p,1}^{1+\frac dp})}\\
\les & \|a^h\|_
{\widetilde L_t^\infty(\dot B_{p,1}^{\frac d p})
\cap
\widetilde L_t^1(\dot B_{p,1}^{\frac d p})}+\|\mathbb P u^h\|_{\widetilde L_t^\infty(\dot B_{p,1}^{-1+\frac dp})\cap\widetilde L_t^1(\dot B_{p,1}^{1+\frac dp})}+\|\mathbb Q u^h\|_{\widetilde L_t^\infty(\dot B_{p,1}^{-1+\frac dp})\cap\widetilde L_t^1(\dot B_{p,1}^{1+\frac dp})}\\
\les& \|a^h\|_
{\widetilde L_t^\infty(\dot B_{p,1}^{\frac d p})
\cap
\widetilde L_t^1(\dot B_{p,1}^{\frac d p})}+\|\mathbb P u^h\|_{\widetilde L_t^\infty(\dot B_{p,1}^{-1+\frac dp})\cap\widetilde L_t^1(\dot B_{p,1}^{1+\frac dp})}+\|w^h\|_{\widetilde L_t^\infty(\dot B_{p,1}^{-1+\frac dp})\cap\widetilde L_t^1(\dot B_{p,1}^{1+\frac dp})}
}
and thus by \eqref{eq3.21} and \eqref{eq3.23}, choosing $j_0$ sufficiently large, we obtain
\EQ{
A^h+U^h\les& \|u_0^h\|_{\dot B_{p,1}^{-1+\frac d p}}
+\|a_0^h\|_{\dot B_{p,1}^{\frac d p}}\\
&+\|(a\cdot\operatorname{div}u)^h\|_{\widetilde L_t^1(\dot B_{p,1}^{\frac d p})}+\|g^h\|_{\widetilde L_t^1(\dot B_{p,1}^{-1+\frac d p})}+\|f^h\|_{\widetilde L_t^1(\dot B_{p,1}^{\frac d p-2})}.\label{eq3.31}
}
It remains to do the nonlinear estimates.

First, we estimate $\|(a\cdot \operatorname{div}u)^h\|_
{\widetilde L_t^1(\dot B_{p,1}^{\frac d p})}$.
We decompose $a\cdot \operatorname{div}u=a \cdot \operatorname{div}u^\ell+a \cdot \operatorname{div}u^h$. By \eqref{eq:prodest1}, we get
\EQ{\label{eq3.29}
\norm{(a \cdot \operatorname{div} u^\ell)^h}_{\tilde{L}_{t}^{1}(\dot B_{p,1}^{\frac{d}{p}})}\les & \norm{a}_{\tilde{L}_t^2(\dot{B}_{p,1}^{\frac{d}{p}})}\norm{u^\ell}_{\tilde{L}_t^2(\dot{B}_{p,1}^{\frac{d}{p}+1})}\les A U \\
\norm{(a \cdot \operatorname{div} u^h)^h}_{\tilde{L}_{t}^{1}(\dot B_{p,1}^{\frac{d}{p}})}\les & \norm{a}_{\tilde{L}_t^\infty(\dot{B}_{p,1}^{\frac{d}{p}})}\norm{u^h}_{\tilde{L}_t^1(\dot{B}_{p,1}^{\frac{d}{p}+1})}\les A U 
}
provided that
\EQ{
2\le q\le p<2d,\quad \theta \leq \frac{d}{q}.\label{eq3.30}
}

Second,  we estimate $\|f^h\|_
{\widetilde L_t^1(\dot B_{p,1}^{\frac d p-2})}$.
Since 
\EQ{\|f^h\|_
{\widetilde L_t^1(\dot B_{p,1}^{\frac d p-2})}\le \|(au)^h\|_
{\widetilde L_t^1(\dot B_{p,1}^{\frac d p-1})}\les \|(au)^h\|_
{\widetilde L_t^1(\dot B_{p,1}^{\frac d p})},
} 
by \eqref{eq:prodest1} we get
\EQ{\|(au)^h\|_
{\widetilde L_t^1(\dot B_{p,1}^{\frac d p})}\les  & \norm{a}_{\tilde{L}_t^2(\dot{B}_{p,1}^{\frac{d}{p}})}\norm{u^\ell}_{\tilde{L}_t^2(\dot{B}_{p,1}^{\frac{d}{p}})}+\norm{a}_{\tilde{L}_t^\infty(\dot{B}_{p,1}^{\frac{d}{p}})}\norm{u^h}_{\tilde{L}_t^1(\dot{B}_{p,1}^{\frac{d}{p}})}\les A U 
}
assuming the condition \eqref{eq3.30}.

Last, we estimate $\|g^h\|_
{\widetilde L_t^1(\dot B_{p,1}^{\frac d p-1})}$.
We write $g(a,u)=-u\cdot \nabla u-R_0(a) a\mathcal{A}u-\nabla (R_1(a)a^2):=g_1+g_2+g_3.$
For $g_1^h$, we have 
\EQ{
\|g_1^h\|_
{\widetilde L_t^1(\dot B_{p,1}^{\frac d p-1})}\les& \|(u\cdot \nabla u)^h\|_
{\widetilde L_t^1(\dot B_{p,1}^{\frac d p-1})}\\
\les& \|(u^\ell\cdot \nabla u)^h\|_
{\widetilde L_t^1(\dot B_{p,1}^{\frac d p-1})}+\|(u^h\cdot \nabla u)^h\|_
{\widetilde L_t^1(\dot B_{p,1}^{\frac d p-1})}.
}
For the term $\|(u^\ell\cdot \nabla u)^h\|_
{\widetilde L_t^1(\dot B_{p,1}^{\frac d p-1})}$, assuming \eqref{eq3.30}, we get by \eqref{eq:prodest2} that
\EQ{
\|(u^\ell\cdot \nabla u)^h\|_
{\widetilde L_t^1(\dot B_{p,1}^{\frac d p-1})}\les U^2.
}
For the term $\|(u^h\cdot \nabla u)^h\|_
{\widetilde L_t^1(\dot B_{p,1}^{\frac d p-1})}$, assuming \eqref{eq3.30}, by \eqref{eq:prodest2} we get
\EQ{
\|(u^h\cdot \nabla u)^h\|_
{\widetilde L_t^1(\dot B_{p,1}^{\frac d p-1})}\les& \norm{u^h}_{\tilde{L}_t^2(\dot{B}_{p,1}^{\frac{d}{p}})}\norm{u^\ell}_{\tilde{L}_t^2(\dot{B}_{p,1}^{\frac{d}{p}})}+\norm{u^h}_{\tilde{L}_t^\infty(\dot{B}_{p,1}^{\frac{d}{p}-1})}\norm{u^h}_{\tilde{L}_t^1(\dot{B}_{p,1}^{\frac{d}{p}+1})}\\
\les& U^2.
}

For $g_2^h$, we have $\|g_2^h\|_
{\widetilde L_t^1(\dot B_{p,1}^{\frac d p-1})}\les \|(R_0(a)a\nabla ^2 u)^h\|_
{\widetilde L_t^1(\dot B_{p,1}^{\frac d p-1})}.$ By the Taylor expansion and \eqref{eq:prodest1} we get
\EQ{
\norm{R_0(a)}_{\dot B_{p,1}^{d/p}}\les 1.
}
We decompose 
$a\nabla ^2 u=a \nabla ^2 u^\ell+a \nabla ^2 u^h$. By \eqref{eq:prodest2}, we get
\EQ{
&\norm{(R_0(a) a \nabla ^2 u^\ell)^h}_{\tilde{L}_{t}^{1}(\dot B_{p,1}^{\frac{d}{p}-1})}\\
&\les \norm{R_0(a)}_{\tilde{L}_{t}^{\infty}(\dot B_{p,1}^{\frac{d}{p}})} \norm{a}_{\tilde{L}_{t}^{2}(\dot B_{p,1}^{\frac{d}{p}})}\norm{u^\ell}_{\tilde{L}_{t}^{2}(\dot B_{p,1}^{\frac{d}{p}+1})}\les AU^\ell
}
and
\EQ{
&\norm{(R_0(a) a \nabla ^2 u^h)^h}_{\tilde{L}_{t}^{1}(\dot B_{p,1}^{\frac{d}{p}-1})}\\
&\les \norm{R_0(a)}_{\tilde{L}_{t}^{\infty}(\dot B_{p,1}^{\frac{d}{p}})} \norm{a}_{\tilde{L}_{t}^{\infty}(\dot B_{p,1}^{\frac{d}{p}})}\norm{u^h}_{\tilde{L}_{t}^{1}(\dot B_{p,1}^{\frac{d}{p}+1})}\les AU^h
}
provided that \eqref{eq3.30} holds.

For $g_3$, we have
\EQ{
\|g_3^h\|_
{\widetilde L_t^1(\dot B_{p,1}^{\frac d p-1})}\les& \|(R_1(a)a^2)^h\|_
{\widetilde L_t^1(\dot B_{p,1}^{\frac d p})}\\
\les& \|a\|_
{\widetilde L_t^2(\dot B_{p,1}^{\frac d p})}\|a\|_
{\widetilde L_t^2(\dot B_{p,1}^{\frac d p})}\les A^2.
}
Then, we obtain
\begin{equation}\label{eq3.42}
\begin{aligned}
\|g^h\|_
{\widetilde L_t^1(\dot B_{p,1}^{\frac d p-1})}\lesssim (A+U)^2.
\end{aligned}
\end{equation}
Combining (\ref{eq3.31}), Case 1-Case~3, and taking $j_0$ sufficiently large such that $2^{-j_0}\ll1$,
we complete the proof of the proposition \ref{prop:high}.
\end{proof}

\subsection{Estimates for $(u^\ell,m^h)$}

In this subsection, we prove the following estimates using the last two equations in \eqref{eq3.1}.

\begin{proposition}\label{prop:ulmh} 
Assume $2 \le p < 2d$. Then we have
\EQ{
U^\ell\les& M^\ell+AU,\\
M^h\les& U^h+AU.
}
\end{proposition}
\begin{proof}
First we estimate $U^\ell:=\norm{u^\ell}_{\widetilde{L}_{t}^\infty({\dot{B}_{p,1}^\frac{d}{p}})\cap \widetilde{L}_{t}^2(\dot{B}_{p,1}^{\frac{d}{p}})}$.

Since $u^\ell =m^\ell-(au)^\ell$, then 
\EQ{
U^\ell \les M^\ell+\norm{(au)^\ell}_{\widetilde{L}_{t}^\infty({\dot{B}_{p,1}^\frac{d}{p}})\cap \widetilde{L}_{t}^2(\dot{B}_{p,1}^{\frac{d}{p}})}
}

We decompose $au=a u^\ell+a u^h$. Then by \eqref{eq:prodest} we have: for $2\leq p<2d$
\EQ{
\norm{(a u^\ell)^{\ell}}_{\widetilde{L}_{t}^\infty({\dot{B}_{p,1}^\frac{d}{p}})\cap \widetilde{L}_{t}^2(\dot{B}_{p,1}^{\frac{d}{p}})}
\les & \norm{a}_{\tilde{L}_t^\infty(\dot{B}_{p,1}^{\frac{d}{p}})}\norm{u^\ell}_{\widetilde{L}_{t}^\infty({\dot{B}_{p,1}^\frac{d}{p}})\cap \widetilde{L}_{t}^2(\dot{B}_{p,1}^{\frac{d}{p}})}\les A U^\ell.
}

\EQ{
\norm{( a u^h)^{\ell}}_{\widetilde{L}_{t}^\infty({\dot{B}_{p,1}^\frac{d}{p}})\cap \widetilde{L}_{t}^2(\dot{B}_{p,1}^{\frac{d}{p}})}\les& \norm{( a u^h)^{\ell}}_{\widetilde{L}_{t}^\infty({\dot{B}_{p,1}^{\frac{d}{p}-1}})}+\norm{( a u^h)^{\ell}}_{\widetilde{L}_{t}^2(\dot{B}_{p,1}^{\frac{d}{p}})}\\
\les& \norm{a}_{\widetilde{L}_{t}^\infty({\dot{B}_{p,1}^{\frac{d}{p}}})}\norm{u^h}_{\widetilde{L}_{t}^\infty({\dot{B}_{p,1}^{\frac{d}{p}-1}})\cap \widetilde{L}_{t}^2(\dot{B}_{p,1}^{\frac{d}{p}})}\les AU^h.
}
Therefore,
\EQ{
U^\ell \les M^\ell + AU.
}

Second, we estimate $M^h:=\norm{m^h}_{
\widetilde L_t^2(\dot B_{p,1}^{\frac{d}{p}})\cap
\widetilde L_t^\infty(\dot B_{p,1}^{\frac{d}{p}-1})}$.

Since $m^h = ((a+1)u)^h = u^h + (au)^h$, then
\EQ{
M^h\les U^h+\norm{(au)^h}_{
\widetilde L_t^2(\dot B_{p,1}^{\frac{d}{p}})\cap
\widetilde L_t^\infty(\dot B_{p,1}^{\frac{d}{p}-1})}.
}
We decompose $au=a u^\ell+a u^h$. Then by \eqref{eq:prodest} we get: for $2\leq p<2d$
\EQ{
\norm{(au^\ell)^h}_{
\widetilde L_t^2(\dot B_{p,1}^{\frac{d}{p}})\cap
\widetilde L_t^\infty(\dot B_{p,1}^{\frac{d}{p}-1})}\les& \norm{au^\ell}_{
\widetilde L_t^2(\dot B_{p,1}^{\frac{d}{p}})\cap
\widetilde L_t^\infty(\dot B_{p,1}^{\frac{d}{p}})}\\
\les& \norm{a}_{
\widetilde L_t^\infty(\dot B_{p,1}^{\frac{d}{p}})}\norm{u^\ell}_{
\widetilde L_t^2(\dot B_{p,1}^{\frac{d}{p}})\cap
\widetilde L_t^\infty(\dot B_{p,1}^{\frac{d}{p}})}\les AU
} 
and
\EQ{
&\norm{(au^h)^h}_{
\widetilde L_t^2(\dot B_{p,1}^{\frac{d}{p}})\cap
\widetilde L_t^\infty(\dot B_{p,1}^{\frac{d}{p}-1})}\\
\les& \norm{au^h}_{
\widetilde L_t^2(\dot B_{p,1}^{\frac{d}{p}})}+\norm{au^h}_{
\widetilde L_t^\infty(\dot B_{p,1}^{\frac{d}{p}-1})}\\
\les& \norm{a}_{
\widetilde L_t^\infty(\dot B_{p,1}^{\frac{d}{p}})}\norm{u^h}_{
\widetilde L_t^2(\dot B_{p,1}^{\frac{d}{p}})}+\norm{a}_{
\widetilde L_t^\infty(\dot B_{p,1}^{\frac{d}{p}})}\norm{u^h}_{
\widetilde L_t^\infty(\dot B_{p,1}^{\frac{d}{p}-1})}\les AU.
} 
Therefore,
\EQ{
M^h\les& U^h+AU.
}
We complete the proof of the proposition. 
\end{proof}

\section{Proof of Theorem \ref{thm}}\label{sect:proofThm}

In this section we prove Theorem~\ref{thm}.  
Assume 
\EQ{\label{eq:allcond}
\frac{2d-1}{q}+\frac{3-d}{2}>0,\ p<2d,\ 2\le q\leq p\leq 2q,\ -\delta(q,q,d)-1<\theta\leq \frac{d}{2q}.
}
which is the intersection of all the conditions in Propositions \ref{prop:low},\ref{prop:high},\ref{prop:ulmh}. 
Suppose that $(a,m,u)$ is a solution to system \eqref{eq3.1} on $[0,T]$
satisfying
\EQ{
\mathcal{X}:=A+M+U\lesssim 1.
}
Let
\EQ{
\mathcal{X}_0=\|(a_0,m_0)^\ell\|_{\dot B_{q,1}^{\frac{2d-1}{q}+\frac{1-d}{2}}\cap \dot B_{q,1}^{\theta+\delta(q,q,d)})}+\|u_0^h\|_{\dot B_{p,1}^{-1+\frac d p}}
+\|a_0^h\|_{\dot B_{p,1}^{\frac d p}}.
}
By Propositions \ref{prop:low},\ref{prop:high},\ref{prop:ulmh}, we obtain
\begin{equation*}
\mathcal{X} \les \mathcal{X}_0 + \mathcal{X}^2.
\end{equation*}
By the standard iteration or bootstrap arguments, assuming $\mathcal{X}_0\ll 1$, we can derive
\begin{equation}\label{xsmall}
\mathcal{X} \leq C\mathcal{X}_{0}.
\end{equation}

Now we take $\theta=\frac{d}{2q}$. By \eqref{eq:semi-qq}, we have $\delta(q,q,d)+1=\frac{d}{q}-\frac{d-1}{2}-$. So we can rewrite the conditions \eqref{eq:allcond} as
\EQ{\label{eq:allcond2}
2\leq q\leq p\leq 2q,\, q<\frac{3d}{d-1}, \, p<2d.
}
Under the above conditions, we have 
\EQ{
\|(a,m)^\ell\|_{\widetilde L_t^\infty(\dot B_{q,1}^{-1+\frac{d}{q}})}
\lesssim&
\|(a_0,m_0)^\ell\|_{\dot B_{q,1}^{\frac{2d-1}{q}-\frac{1+d}{2}}}+\|h(a,m)^\ell\|_{\widetilde L_t^1(\dot B_{r,1}^{\frac{2d-1}{q}+\frac{-1-d}{2}})}\\
\les& \|(a_0,m_0)^\ell\|_{\dot B_{q,1}^{\frac{2d-1}{q}-\frac{1+d}{2}}}+\|h(a,m)^\ell\|_{\widetilde L_t^1(\dot B_{r,1}^{\theta+\delta(q,q,d)})}\\
\les & \|(a_0,m_0)^\ell\|_{\dot B_{q,1}^{\frac{2d-1}{q}-\frac{1+d}{2}}}+AM.
}
Using the last equation of \eqref{eq3.1}, we get
\EQ{
\|u^\ell\|_{\widetilde L_t^\infty(\dot B_{p,1}^{-1+\frac{d}{p}})}\les \|m^\ell\|_{\widetilde L_t^\infty(\dot B_{q,1}^{-1+\frac{d}{q}})}+\|(au)^\ell\|_{\widetilde L_t^\infty(\dot B_{p,1}^{-1+\frac{d}{p}})}\les M+A(\|u^\ell\|_{\widetilde L_t^\infty(\dot B_{p,1}^{-1+\frac{d}{p}})}+U)
}
which implies 
\EQ{
\|u^\ell\|_{\widetilde L_t^\infty(\dot B_{p,1}^{-1+\frac{d}{p}})}\ll 1
}
if $\|(a_0,m_0)^\ell\|_{\dot B_{q,1}^{\frac{2d-1}{q}-\frac{1+d}{2}}}\ll 1$.
Therefore, we have shown $\norm{(a,u)}_{C_T\X_p}\ll 1$, and the solution $(a,u)$ exists globally.

\section*{Acknowledgments}
No generative artificial intelligence (AI) was used in the writing or preparation of this manuscript. R. Chen gratefully acknowledges the financial support of a Monash Graduate Scholarship (MGS) and FT230100588. Z. Guo is supported by ARC FT230100588 and DP260100485. 






\end{document}